\documentclass{article}

\usepackage{amsmath, amsthm, amsfonts}

\usepackage{mdframed}
\theoremstyle{definition}
\newmdtheoremenv{boxProb}{Problem}
\newmdtheoremenv{boxDef}{Definition}
\newmdtheoremenv{boxCor}{Corollary}
\newmdtheoremenv{boxThm}{Theorem}
\newmdtheoremenv{compjob}{Computational Job}
\newmdtheoremenv{reqi}{Requirement}

\usepackage{graphicx}
\usepackage{subfigure}          

\usepackage{placeins}

\usepackage{tabularx}

\usepackage{lmodern,textcomp}

\usepackage{algorithm}
\usepackage{algpseudocode}

\usepackage{enumitem}

\usepackage{xspace} 
\newcommand\largeparbreak{\par\bigskip}

\newcommand*\tageq{\refstepcounter{equation}\tag{\theequation}}

\newcommand{\blambda}{\boldsymbol{\lambda}\xspace}
\newcommand{\bxi}{\boldsymbol{\xi}\xspace}

\renewcommand{\t}{^\textsf{T}\xspace}

\newcommand{\away}[1]{}

\newcommand{\R}{\mathbb{R}\xspace}

\newcommand{\N}{\mathbb{N}\xspace}

\newcommand{\cL}{\mathcal{L}\xspace}

\newcommand{\cO}{\mathcal{O}\xspace}

\newcommand{\bB}{\textbf{B}\xspace}

\newcommand{\bS}{\textbf{S}\xspace}

\newcommand{\bK}{\textbf{K}\xspace}

\newcommand{\bx}{\textbf{x}\xspace}

\newcommand{\by}{\textbf{y}\xspace}

\newcommand{\bei}[1]{{\textbf{e}}\xspace}
\newcommand{\bw}{\textbf{w}\xspace}

\newcommand{\bI}{\textbf{I}\xspace}

\newcommand{\bO}{\textbf{0}\xspace}

\newcommand{\tblambda}{\tilde{\boldsymbol{\lambda}}\xspace}

\newcommand{\tomega}{{\tilde{\omega}}\xspace}

\newcommand{\tol}{{\textsf{tol}}\xspace}

\title{Modified Augmented Lagrangian Method for the minimization of functions with quadratic penalty terms}

\author{Martin Neuenhofen}

\date{\today}

\begin{document}

\maketitle

\begin{abstract}
	The Augmented Lagrangian Method (ALM) is an iterative method for the solution of equality-constrained non-linear programming problems. In contrast to the quadratic penalty method, the ALM can satisfy equality constraints in an exact way. Further, ALM is claimed to converge in less iterations, indicating that it is superior in approach to a quadratic penalty method.
	
	It is referred to as an advantage that the ALM solves equality constraints in an exact way, meaning that the penalty parameter does not need to go to infinity in order to yield accurate feasibility for the constraints. However, we sometimes actually want to minimize an unconstrained problem that has large quadratic penalty terms. In these situations it is unclear how the ALM could be utilized in the correct way. This paper presents the answer: We derive a modified version of the ALM that is also suitable for solving functions with large quadratic penalty terms.
\end{abstract}


\section{Introduction}

We consider numerical methods for the minimization of three types of optimization problems. These problems are stated below.

\paragraph{Merit problem}
We consider the unconstrained minimization of a function with quadratic penalty terms:
\begin{align*}
	\min_{\bx \in \R^n} \quad \Phi(\bx) := f(\bx) + \frac{1}{2 \cdot \omega} \cdot \|c(\bx)\|_2^2\tageq\label{eqn:MeritProblem}\,,
\end{align*}
where $m,n \in \N$, a very small number $\omega \in \R^+\setminus\lbrace 0 \rbrace$, $f:\ \R^n \rightarrow \R$ and $c:\ \R^n \rightarrow \R^m$ are given. Notice that the statement of this problem makes sense regardless whether $m\leq n$ or $m>n$\,. This is an unconstrained optimization problem. Under suitable assumptions on $\Phi$ one can show well-posedness.

\paragraph{Penalty problem}
Though redundant, it helps the presentation to define also the following optimization problem:
\begin{align*}
	\min_{\bx \in \R^n} \quad \Psi(\bx) := f(\bx) + \frac{1}{2 \cdot \tomega} \cdot \|c(\bx)\|_2^2\tageq\label{eqn:PenaltyProblem}
\end{align*}
where $\tomega \in \R^+\setminus\lbrace 0 \rbrace$ is a small number. We think of $\tomega$ as a number that is small but that is still larger than $\omega$. E.g., it could be $\omega = 10^{-10}$ and $\tomega = 10^{-2}$\,. Given a suitable initial guess, penalty problems can be considered solvable with a conventional Newton-type method whereas merit problems can be considered intractable because $\omega$ is usually so small that their scaling is too bad.

\paragraph{Constrained problem}
We consider the following constrained problem.
\begin{subequations}
	\begin{align*}
		\min_{\bx \in \R^n}& & 		 f(\bx)& \tageq 	\\
		\text{s.t.}& 		& 		 c(\bx)& = \bO\tageq
	\end{align*}\label{eqn:ConstrainedProblem}
\end{subequations}
If a solution $\bx$ of $c(\bx)=\bO$ exists then again under suitable assumptions one can show well-posedness. Further, in this case solutions of problem \eqref{eqn:MeritProblem} converge to solutions of \eqref{eqn:ConstrainedProblem} for decreasing positive values of $\omega$\,. Since in practice we only solve optimization problems up to some tolerance $\tol>0$ for the residual of the Karush-Kuhn-Tucker (KKT) conditions, we can consider \eqref{eqn:ConstrainedProblem} as a sub-class of problem \eqref{eqn:MeritProblem}. We do so because when choosing $\omega \in \cO(\tol)$, a solution of \eqref{eqn:MeritProblem} is also a sufficiently accurate solution of \eqref{eqn:ConstrainedProblem}, cf. \cite{LOQO}.

The Lagrangian of \eqref{eqn:ConstrainedProblem} is
\begin{align}
	\cL(\bx,\blambda) := f(\bx) - \blambda\t\cdot c(\bx)\,.\label{eqn:Lagrangian}
\end{align}

\paragraph{Outline}
In this paper we consider three local solution methods:
\begin{enumerate}
	\item The penalty method: It is an iterative method for solving \eqref{eqn:MeritProblem} by solving a sequence of \eqref{eqn:PenaltyProblem} with decreasing values for $\tomega$. Since from a numerical perspective \eqref{eqn:ConstrainedProblem} is a special case of \eqref{eqn:MeritProblem}, the penalty method can also be used to solve \eqref{eqn:ConstrainedProblem}.
	\item The ALM: It is a method for solving \eqref{eqn:ConstrainedProblem} by solving a finite sequence of \eqref{eqn:PenaltyProblem} with a value for $\tomega$ that remains moderate, e.g. $10^{-2}$, throughout the sequence.
	\item A modified ALM: It is a method for solving \eqref{eqn:MeritProblem} by solving a finite sequence of \eqref{eqn:PenaltyProblem} with a value for $\tomega$ that remains small, e.g. $10^{-2}$, throughout the sequence.
\end{enumerate}

The first and second method are well-known from the literature. The third approach is a novelty that we propose in this paper. The motivation for using this method is that it can solve \eqref{eqn:MeritProblem} by solving a sequence of \eqref{eqn:PenaltyProblem}. But in contrast to the penalty method it uses moderate values of $\tomega$, as is the case of the conventional ALM. This is beneficial, as it is claimed that ALM usually converges in less iterations than the penalty method \cite{ALM}.

\section{Summary of penalty method and ALM}
We review the two methods from the literature. We work out their computational framework and give pseudo-code that can be used for a practical implementation.

\subsection{Penalty method}
The penalty method solves a sequence of \eqref{eqn:PenaltyProblem}, where $\tomega$ iteratively decreases until it finally reaches $\omega$\,. The solution of each problem \eqref{eqn:PenaltyProblem} is used as initial guess to solve the subsequent problem \eqref{eqn:PenaltyProblem} with a decreased value of $\tomega$\,. While decreasing values of $\tomega$ make the problem \eqref{eqn:PenaltyProblem} more difficult to solve, the initial guesses obtained from the former solution yield a compensating effect. In particular, the found minimizers become more and more accurate initial guesses for the subsequent penalty problem, altogether resulting in a potentially efficient approach for solving \eqref{eqn:MeritProblem}.

The method can be described with the following pseudocode.

\begin{algorithmic}[1]
\State Given is: $f$, $c$, $\omega>0$, $0<\theta_\omega<1$, $\tomega_1>\omega$, $\bx_0 \in \R^n$
\For{$k=1,2,3,...$}
	\State Solve \eqref{eqn:PenaltyProblem} for $\tomega=\tomega_k$ with initial guess $\bx_{k-1}$.
	\State Write the solution of \eqref{eqn:PenaltyProblem} into the vector $\bx_{k} \in \R^n$\,.
	\If{$\tomega_k$==$\omega$}
		\State \Return $\bx_{k}$
	\Else
		\State Update $\tomega_{k+1} := \max\lbrace\,\theta_\omega \cdot \tomega_k\,,\,\omega\,\rbrace$\,.
	\EndIf
\EndFor
\end{algorithmic}

\paragraph{The Newton iteration within the penalty method}
In order to solve \eqref{eqn:PenaltyProblem}, typically a Newton-type method is used. In the simplest approach, the Newton type method needs two things:
\begin{itemize}
	\item Root function: This is a multivariate function that has as many output as input dimensions. Its root, computed by Newton's method, determines a local minimizer.
	\item A merit-function: This is a function for which the Newton step (subject to suitable modifications) is a descent direction. It is needed within the line search, which in turn is needed for globalization of the Newton method.
\end{itemize}

Problem \eqref{eqn:PenaltyProblem} has the Karush-Kuhn Tucker (KKT) \cite{NumOpt} condition
\begin{align*}
	\nabla\Psi(\bx) \equiv \nabla f(\bx) - \nabla c(\bx) \cdot \left(\frac{-1}{\tomega} \cdot c(\bx)\right) = \bO\,.
\end{align*}
As discussed in \cite{Amudson}, for sake of well-posedness it is advantageous to rather consider the following equivalent set of equations, where $\tblambda=-1/\tomega \cdot c(\bx)$ is a substituted auxiliary variable.
\begin{align*}
	F(\bx,\tblambda) := \begin{pmatrix}
	\nabla f(\bx) - \nabla c(\bx) \cdot \tblambda \\
	c(\bx) + \tomega \cdot \tblambda
	\end{pmatrix} = \bO
\end{align*}
The Jacobian is
\begin{align*}
	DF(\bx,\tblambda) = \begin{bmatrix}
		\nabla_{\bx\bx}^2 \cL(\bx,\tblambda) &-\nabla c(\bx)\\
		\nabla c(\bx)\t & \tomega \cdot \bI
	\end{bmatrix}\,.
\end{align*}
A symmetric indefinite system is obtained by multiplying the second column with $-1$\,. Further, replacing $\nabla_{\bx\bx}^2 \cL(\bx,\tblambda)$ with a symmetric matrix $\bB \in \R^{n \times n}$ we obtain the symmetric Newton system
\begin{align}
	\underbrace{\begin{bmatrix}
	\bB 			&\nabla c(\bx)\\
	\nabla c(\bx)\t &-\tomega \cdot \bI
	\end{bmatrix}}_{=:\bK} \cdot \begin{pmatrix}
	\Delta\bx\\
	-\Delta\tblambda
	\end{pmatrix} = -\begin{pmatrix}
	\nabla_\bx \cL(\bx,\tblambda)\\
	c(\bx) + \tomega \cdot \tblambda
	\end{pmatrix}\,.\label{eqn:NewtonSystem_PenaltyProblem}
\end{align}
We use the short-hands $\bw := (\bx\t,\tblambda\t)\t$, $\Delta\bw := (\Delta\bx\t,\Delta\tblambda\t)\t \in \R^{n+m}$\,. The symmetric Newton matrix we call $\bK$.

We say that $\bB$ is chosen \textit{suitable} when $\bK$ has $n$ strictly positive and $m$ strictly negative eigenvalues. For instance, $\bB$ is suitable when it is positive definite. This can be achieved in the simplest way by adding a multiple of the identity to $\nabla^2_{\bx\bx}\cL(\bx,\tblambda)$. More sophisticated is the use of a symmetric indefinite LDL-factorization with a subsequent modification of the diagonal matrix. For details we refer to \cite[Sec.~4.1]{ForsgrenGill}. We assume in the following that $\bB$ is suitable.

Finally, we need the merit function. As shown in \cite{ForsgrenGill}, $\Delta\bw$ is a descent direction for the following merit function for arbitrary values of $\nu \in \R^+\setminus\lbrace 0 \rbrace$\,:
\begin{align*}
	M(\bw) := f(\bx) + \frac{\tomega}{2} \cdot \|c(\bx)\|_2^2 + \frac{\nu}{2\cdot\tomega} \cdot \|c(\bx) + \tomega \cdot \tblambda\|_2^2\,. \tageq\label{eqn:MeritPenaltyMethod}
\end{align*}

For the purpose of an overview and a concrete method, below we state the solution procedure for \eqref{eqn:PenaltyProblem}, inherited in line 2 of the penalty method, below.

\begin{algorithmic}[1]
	\State Given is: $f,c$, $\tomega>0$, initial guess $\bx$, $\tol>0$\,.
	\State $\tblambda := -1/\tomega \cdot c(\bx)$
	\While{$\|F(\bw)\|_2 >\tol$}
		\State Choose $\bB$ suitable; e.g. $\bB := \nabla_{\bx\bx}^2 \cL(\bx,\tblambda) + \xi \cdot \bI \succ \bO$ with $\xi\geq 0$ sufficiently large.
		\State Solve the linear system \eqref{eqn:NewtonSystem_PenaltyProblem}.
		\State Find $\alpha \in (0,1]$ that satisfies a sufficient descent condition for the merit function in \eqref{eqn:MeritPenaltyMethod}.
		\State Update $\bw := \bw + \alpha \cdot \Delta\bw$\,.
	\EndWhile
\end{algorithmic}

\subsection{ALM}\label{Sec:ALM}
The ALM solves a sequence of \eqref{eqn:PenaltyProblem}, where usually $\tomega$ remains bounded below by some moderate number. The actual requirements on the minimum value for $\tomega$ remain technical and are discussed to in \cite{ALM} and the references therein.

In order for minimizers of \eqref{eqn:PenaltyProblem} to solve \eqref{eqn:ConstrainedProblem}, the function $f$ in $\Psi$ is replaced by $\cL(\bx,\blambda)$, where $\blambda \in \R^m$ is computed suitable through an outer iteration scheme. For better reference we introduce an outer iteration index $k \in \N_0$ that we add as footnote to $\blambda_k$\,.

For the $k$-th outer iteration we write out the particular problem of kind \eqref{eqn:PenaltyProblem} that is solved within ALM:
\begin{align*}
\min_{\bx \in \R^n} \quad \Psi_k(\bx) := \cL(\bx,\blambda_k) + \frac{1}{2\cdot\tomega} \cdot \|c(\bx)\|_2^2\tageq\label{eqn:ALM_PenaltyProblem}
\end{align*}
In this problem, $\blambda_k$ is a fixed given vector from the outer iteration $k$. $\tomega$ is assumed to be sufficiently small so that the problem is bounded.

As we did for the penalty problem solved within the penalty method, from the solution of \eqref{eqn:ALM_PenaltyProblem} we obtain a tuple $(\bx,\tblambda) \in \R^n \times \R^m$, where $\tblambda \approx -1/\tomega \cdot c(\bx)$, where $\approx$ accounts for the inaccuracy in terms of the tolerance $\tol>0$\,.

The algorithm of ALM goes as follows.
\begin{algorithmic}[1]
	\State Given is: $f$, $c$, $0<\theta_\tomega<1$, $0<\theta_\lambda<1$, $\tomega_1>0$, $\bx_0 \in \R^n$
	\State $\tblambda_0 := -1/\tomega_1 \cdot c(\bx_0)$\,,\quad $\blambda_0 := \bO$
	\For{$k=1,2,3,...$}
		\State Solve \eqref{eqn:ALM_PenaltyProblem} for $\tomega=\tomega_k$ with initial guess $(\bx_{k-1},\tblambda_{k-1})$.
		\State Write the solution of \eqref{eqn:ALM_PenaltyProblem} into the tuple $(\bx_{k},\tblambda_k) \in \R^n \times \R^m$\,.
		\If{$\|\nabla_{\bx\blambda}\cL(\bx_k,\blambda_k)\|_2\leq \tol$}
			\State \Return $\bx_{k}$
		\Else
			\State // Update $\blambda$ and $\tomega$
			\If{$\|c(\bx_k)\|_2 \leq \theta_\lambda \cdot \min_{j=0,...,k-1}\lbrace\,\|c(\bx_j)\|_2\,\rbrace$}
				\State $\blambda_{k+1} := \blambda_k + \tblambda_k$\,,\quad $\tomega_{k+1} := \tomega_k/\sqrt{\theta_\omega}$
				\State $\tblambda_k := \bO$
			\Else
				\State $\blambda_{k+1} := \blambda_k$\,,\quad $\tomega_{k+1} = \theta_\omega \cdot \tomega_k$
			\EndIf
		\EndIf
	\EndFor
\end{algorithmic}
This algorithm, apart from constant parameters $\theta_\omega,\theta_\lambda$, e.g. $\theta_\omega=\theta_\lambda=0.1$\,, is equivalent to the ALM as presented in \cite{ALM}. The increase of $\tomega$ in the update in line 11 is optional.

We describe the steps of the algorithm. We start by solving problem \eqref{eqn:ALM_PenaltyProblem}. For moderate but sufficiently small values of $\tomega$ the solution $\bx_k$ of this problem will be bounded and yields a small norm for $\|c(\bx_k)\|_2$\,. In line 6 the KKT conditions are checked. If they are satisfied up to a tolerance $\tol>0$ then the problem is considered solved and the solution is returned. Otherwise, $\blambda$ and $\tomega$ are updated. This can happen in either of two ways, of which we first consider the update in line 11. ALM considers $\tblambda_k$ as a guess for an update of the Lagrange multiplicator $\blambda$. This can be motivated from similarities of the KKT conditions for the four problems \eqref{eqn:MeritPenaltyMethod},\eqref{eqn:PenaltyProblem},\eqref{eqn:ConstrainedProblem} and \eqref{eqn:ALM_PenaltyProblem}. The update for $\blambda$ in line 11 however requires that $\tblambda$ converges to zero for increasing indices $k$; clearly, otherwise $\blambda$ could never converge. For sufficiently small (but bounded below) values of $\tomega$ this will eventually happen under mild assumptions, cf. \cite{ALM} and references therein. However, during earlier iterations the value of $\tomega$ may be too small. Thus there is a conditional statement. If the norm of $\|c(\bx_k)\|_2$ decreases monotonously then the whole iteration is considered convergent, i.e. $\tomega$ is assumed to be sufficiently large. In this case the ALM update is chosen. Otherwise the iteration resorts to a penalty method (update in line 14), which converges under milder assumptions -- until $\|c(\bx_k)\|_2$ decreases again.

The reason why $\tblambda_k$ will eventually converge to zero is due to the values of  $\blambda_k$. In ALM $\blambda_k$ is an iteratively improving guess for the solution of the local Lagrange-dual problem \cite{Boyd}. $\blambda_k$ balances the constraints against the objective. Therefore, eventually for large indices of the outer iteration $k$ the local minimizer $\bx_k$ of
\begin{align*}
	f(\bx_k) - \blambda_k\t\cdot c(\bx_k)
\end{align*}
will just satisfy $c(\bx_k)=\bO$. It is identical to a local minimizer of \eqref{eqn:ALM_PenaltyProblem}, regardless whether $\tomega$ is small or large (which is why one can even increase it). Hence, $\tblambda_k = -1/\tomega \cdot c(\bx_k)$ converges to zero.

\paragraph{Newton system and merit function in the ALM}
As discussed formerly for the penalty method, in order to solve the penalty problem \eqref{eqn:ALM_PenaltyProblem}, we typically use a Newton method. This requires the following two items:
\begin{itemize}
	\item Root function
	\item Merit-function
\end{itemize}

Substituting again $\tblambda = -\tomega \cdot c(\bx)$, the KKT conditions of \eqref{eqn:ALM_PenaltyProblem} can be expressed equivalently as
\begin{align*}
	F(\bx,\tblambda) := 
	\begin{pmatrix}
		\nabla f(\bx) - \nabla c(\bx) \cdot (\blambda_k + \tblambda) \\
		c(\bx) + \tomega \cdot \tblambda
	\end{pmatrix} = \bO\,.\tageq\label{eqn:ALM_res}
\end{align*}
Notice that $\blambda_k$, the given fixed Lagrange multiplier estimate from the outer iteration $k$, appears as a constant in the root function. Using again a short-hand $\bw$ for the argument of $F$, a well-posed Newton system for the update $\Delta\bw = (\Delta\bx\t,\Delta\blambda\t)\t$ is given by
\begin{align*}
	\underbrace{\begin{bmatrix}
		\bB & \nabla c(\bx)\\
		\nabla c(\bx)\t & -\tomega \cdot \bI
	\end{bmatrix}}_{=:\bK} \cdot \begin{pmatrix}
		\Delta\bx\\
		-\Delta\tblambda
	\end{pmatrix}
	= -\begin{pmatrix}
		\nabla_\bx\cL(\bx,\blambda_k + \tblambda)\\
		c(\bx) + \tomega \cdot \tblambda
	\end{pmatrix}\,,
\end{align*}
where this time $\bB$ must be a suitable approximation of $\nabla_{\bx\bx}^2\cL(\bx,\blambda_k+\tblambda)$\,. Again, with suitable we mean that $\bK$ must have $n$ positive and $m$ negative eigenvalues.

The step $\Delta\bw$ is a descent direction for the merit function
\begin{align*}
	M_k(\bw) = f(\bx) - \blambda_k\t\cdot c(\bx) + \frac{1}{2\cdot \tomega} \cdot \|c(\bx)\|^2_2 + \frac{\nu }{2\cdot \tomega}\cdot \|c(\bx) + {\tomega} \cdot \tblambda\|_2^2
\end{align*}
for a fixed arbitrary $\nu>0$ \cite[Sec.~3.1]{ForsgrenGill}. The algorithm for the Newton iteration within the ALM is then analogous to the formerly discussed Newton iteration within the penalty method: Given $\bw$, we first solve the Newton system for $\Delta\bw$ and then perform a line search on $M_k$. We repeat these two steps until $\|F(\bw)\|_2 \leq \tol$ is satisfied for $F$ from \eqref{eqn:ALM_res}.

\section{Modified Augmented Lagrangian Method for the penalty problem}

\paragraph{Motivation}
Both the penalty method and the ALM can be used to solve the constrained problem \eqref{eqn:ConstrainedProblem}. In reaching for this purpose, the ALM is claimed to be more efficient than the penalty method \cite{ALM}. This is for several arguments:
\begin{itemize}
	\item In the penalty method the parameter $\tomega$ must be chosen very small whereas in the ALM it often remains at a moderate size.
	\item A very small value of $\tomega$ can result in numerical problems when solving the Newton system. Thus, penalty methods are considered less numerically robust than ALM.
	\item A very small value of $\tomega$ can also result in steep valleys in the surface of the merit function. Especially for non-linear functions $c$ this valley is curved, like the Rosenbrock "banana" function. This leads to the problem that the line search will make slower progress in the penalty method.
\end{itemize}

On the other hand, the penalty method can solve the more general problem \eqref{eqn:MeritProblem}. We would like to adapt ALM such that it can be used for solving this problem, too. To this end we derive a modified ALM that is suitable for solving both problems \eqref{eqn:MeritProblem} and \eqref{eqn:ConstrainedProblem}.

Our modified ALM works with a method parameter $\omega \in \R^+ \cap \lbrace 0 \rbrace$. It is a true generalization of the conventional ALM as presented above. This is meant in the following sense: For $\omega>0$ the modified ALM solves the problem \eqref{eqn:MeritProblem}. If $\omega=0$ then the modified ALM becomes equivalent to the conventional ALM and solves \eqref{eqn:ConstrainedProblem}. Further, the modification is only a small change that can be easily incorporated into existing implementations of ALM.

\subsection{Derivation of modified ALM}
To derive the modified ALM, we proceed in two steps:
\begin{enumerate}
	\item We write the conventional ALM as a template.
	\item We then fit the merit problem into this template and simplify the resulting expressions.
\end{enumerate} 

\paragraph{The template of ALM}
The ALM can be applied to programming problems that can be expressed as the constrained problem
\begin{subequations}
\begin{align}
	\min_{\by \in \R^N}& 	&  	\tilde{f}(\by)& \\
	\text{s.t.}& 			&	\tilde{c}(\by)&=\bO \,.
\end{align}\label{eqn:TemplateContrainedProblem}
\end{subequations}
Then the root-function of ALM is
\begin{align*}
	F(\by,\tblambda) := 
		\begin{pmatrix}
			\nabla_\by \tilde{f}(\by) - \nabla_\by \tilde{c}(\by) \cdot (\blambda_k + \tblambda) \\
			\tilde{c}(\by) + \tomega \cdot \tblambda
		\end{pmatrix} = \bO\,.
\end{align*}
We use the iterate $\bw := (\by\t,\tblambda\t)\t$. The merit function is
\begin{align*}
	M_k(\bw) = \tilde{f}(\by) - \blambda_k\t\cdot \tilde{c}(\by) + \frac{1}{2\cdot\tomega} \cdot \|\tilde{c}(\by)\|^2_2 + \frac{\nu}{2\cdot\tomega}\cdot \|\tilde{c}(\by) + \tomega \cdot \tblambda\|_2^2
\end{align*}
for $\bw:=(\by,\tblambda)$. And the update formulas, as in the ALM algorithm lines 11 and 14, are:
\begin{align*}
	\blambda_{k+1} &:= \blambda_k + \tblambda_k\,, & \tomega_{k+1} &:= \tomega_k / \sqrt{\theta_\omega}\quad& &\text{if }\|\tilde{c}(\by_k)\|_2 \leq \theta_\lambda \cdot \min_{0\leq j < k}\lbrace\,\|\tilde{c}(\by_j)\|_2\,\rbrace\\
	\blambda_{k+1} &:= \blambda_k\,, & \tomega_{k+1} &:= {\theta_\omega} \cdot \tomega_k\quad\quad & &\text{otherwise }
\end{align*}

\paragraph{Expressing Merit-Problem as Constrained Problem}
We can express \eqref{eqn:MeritProblem} equivalently as a constrained problem by using auxiliary variables $\bxi = 1/\omega \cdot c(\bx)$:
\begin{subequations}
	\begin{align}
		\min_{\bx \in \R^n, \bxi \in \R^m}& 	& 	&f(\bx) + \frac{\omega}{2} \cdot \|\bxi\|_2^2\\
		\text{s.t.}& 							& 	&c(\bx) + \omega \cdot \bxi = \bO
	\end{align}\label{eqn:UsedTemplateCP}
\end{subequations}
This poblem in turn fits into the template \eqref{eqn:TemplateContrainedProblem} by using the following definitions:
\begin{align*}
	\by &:= (\bx\t,\bxi\t)\t \in \R^N\,,\quad N:= n+m\\
	\tilde{f}(\by) &:= f(\bx) + \frac{\omega}{2} \cdot \|\bxi\|_2^2\\
	\tilde{c}(\by) &:= c(\bx) + \omega \cdot \bxi
\end{align*}

From the template we obtain the root function
\begin{align}
	F(\bx,\bxi,\tblambda) := 
	\begin{pmatrix}
		\begin{pmatrix}
			\nabla f(\bx)\\
			\omega \cdot \bxi
		\end{pmatrix} - \begin{pmatrix}
			\nabla c(\bx)\\
			\omega \cdot \bI
		\end{pmatrix} \cdot (\blambda_k + \tblambda) \\
		c(\bx) + \omega \cdot \bxi + \tomega \cdot \tblambda
	\end{pmatrix} = \bO\label{eqn:F}
\end{align}
and the merit function
\begin{align*}
	M_k(\bx,\bxi,\tblambda) = &f(\bx) + \frac{\omega}{2} \cdot \|\bxi\|_2^2 -\blambda_k\t \cdot \big(\,c(\bx) + \omega \cdot \bxi\,\big)\\
	&\quad+ \frac{1}{2\cdot\tomega} \cdot \|c(\bx) + \omega \cdot \bxi\|_2^2 + \frac{\nu}{2\cdot\tomega} \cdot \|c(\bx) + \omega \cdot \bxi + \tomega \cdot \tblambda\|_2^2\,.\tageq\label{eqn:M}
\end{align*}

\paragraph{Elimination of $\bxi$}
We propose to eliminate $\bxi$. From the second row in $F$ in \eqref{eqn:F} we find $\bxi = \blambda_k + \tblambda$. Inserting this for \eqref{eqn:F}, we arrive at
\begin{align}
	F(\bx,\tblambda) := 
	\begin{pmatrix}
		\nabla f(\bx)-\nabla c(\bx) \cdot (\blambda_k + \tblambda) \\
		c(\bx) + \omega \cdot \blambda_k + \left(\omega+\tomega\right) \cdot \tblambda
	\end{pmatrix} = \bO\,.\label{eqn:Newton_residual}
\end{align}
The modified symmetric Newton system is then
\begin{align}
	\underbrace{\begin{bmatrix}
		\bB & \nabla c(\bx)\\
		\nabla c(\bx)\t & -(\omega+\tomega) \cdot \bI
	\end{bmatrix}}_{=:\bK} \cdot 
	\begin{pmatrix}
		\Delta\bx\\
		-\Delta\tblambda
	\end{pmatrix} =
	-\begin{pmatrix}
		\nabla_\bx \cL(\bx,\blambda_k+\tblambda)\\
		c(\bx) + \omega \cdot \blambda_k + (\omega + \tomega) \cdot \tblambda
	\end{pmatrix}\,,\label{eqn:Newton_red}
\end{align}
where $\bB$ is a suitable approximation of $\nabla^2_{\bx\bx}\cL(\bx,\blambda_k+\tblambda)$\,, i.e. so that $\bK$ has $n$ positive and $m$ negative eigenvalues.

We also insert the elimination of $\bxi$ into the merit function \eqref{eqn:M}:
\begin{align*}
	M_k(\bx,\tblambda) 
= 	&f(\bx) 
	+ \frac{\omega}{2} \cdot \|\blambda_k + \tblambda\|_2^2
	- \blambda_k\t \cdot \big(\, c(\bx) + \omega \cdot \blambda_k + \omega \cdot \tblambda\,\big)\\
	& \quad \frac{1}{2\cdot\tomega} \cdot \|c(\bx) + \omega \cdot \blambda_k + \omega \cdot \tblambda\|_2^2\\
	&\quad + \frac{\nu}{2 \cdot\tomega} \cdot \|c(\bx) + \omega \cdot \blambda_k + (\omega + \tomega) \cdot \tblambda\|_2^2
\end{align*}
We can use the following identity:
\begin{align*}
	 &\frac{\omega}{2} \cdot \|\blambda_k + \tblambda\|_2^2 - \blambda\t\cdot\big(\,c(\bx) + \underbrace{\omega \cdot \blambda}_{constant} + \omega \cdot \tblambda\,\big)\\
	=&\frac{\omega}{2} \cdot \|\tblambda\|_2 + \omega \cdot \blambda_k\t \cdot \tblambda - \blambda_k\t \cdot c(\bx) - \omega \cdot \blambda_k\t \cdot \tblambda + const \\
	=&\frac{\omega}{2} \cdot \|\tblambda\|_2^2 - \blambda_k\t\cdot c(\bx) + const
\end{align*}
Using the identity, the merit function can be simplified to
\begin{align*}
	M_k(\bx,\tblambda) :=& f(\bx) + \frac{\omega}{2} \cdot \|\tblambda\|_2^2 - \blambda_k\t \cdot c(\bx) + \frac{1}{2\cdot\tomega} \cdot  \|c(\bx) + \omega \cdot \blambda_k + \omega \cdot \tblambda\|_2^2\\
	&\quad + \frac{\nu }{2\cdot\tomega} \cdot \|c(\bx)+ \omega \cdot \blambda_k + (\omega + \tomega) \cdot \tblambda\|_2^2\,. \tageq\label{eqn:M_red}
\end{align*}
We observe the desirable property that $\omega>0$ forces boundedness for $\tblambda$.

\paragraph{Convergence of the modified ALM}
The modified ALM inherits all convergence properties of ALM. This is because the modified ALM is algebraically equivalent to applying the ALM to problem \eqref{eqn:UsedTemplateCP}.

Modified ALM can be considered a true generalization of ALM because modified ALM solves \eqref{eqn:Newton_residual} where $\blambda_k$ is chosen such that in the limit $\tblambda = \bO$. ALM on the other hand solves \eqref{eqn:ALM_res} where $\blambda_k$ is chosen such that in the limit $\tblambda = \bO$. For $\omega=0$ the modified ALM becomes equivalent to ALM. For $\omega>0$ it can also solve the problem \eqref{eqn:MeritProblem}, which as formerly discussed is a generalization of \eqref{eqn:ConstrainedProblem}.

\subsection{Algorithm of modified ALM}
In summary we state the algorithm of our proposed modified ALM:
\begin{algorithmic}[1]
	\State Given is: $f$, $c$, $\omega>0$, $0<\theta_\omega<1$, $0<\theta_\lambda<1$, $\tomega_1>\omega$, $\bx_0 \in \R^n$
	\State $\tblambda_0 := -1/\tomega_1 \cdot c(\bx_0)$\,,\quad $\blambda_0 := \bO$
	\For{$k=1,2,3,...$}
		\State Solve \eqref{eqn:Newton_residual} for $\tomega=\tomega_k$ with initial guess $(\bx_{k-1},\tblambda_{k-1})$, using the Newton iteration defined by system \eqref{eqn:Newton_red} and merit \eqref{eqn:M_red}.
		\State Write the solution of \eqref{eqn:Newton_residual} into the tuple $(\bx_{k},\tblambda_k) \in \R^n \times \R^m$\,.
		\If{$\max\lbrace\,\omega +\tomega_k\,,\,\|\nabla c(\bx_k)\|_2\,\rbrace \cdot \|\tblambda\|_2\leq \tol$}
			\State \Return $\bx_{k}$
		\Else
			\State // Update $\blambda$ and $\tomega$
			\If{$\|c(\bx_k) + \omega \cdot \blambda_k + \omega \cdot \tblambda_k\|_2 \leq \theta_\lambda \cdot \min_{0\leq j<k}\lbrace\,\|c(\bx_j) + \omega \cdot \blambda_j + \omega \cdot \tblambda_j\|_2\,\rbrace$}
				\State $\blambda_{k+1} := \blambda_k + \tblambda_k$\,,\quad $\tomega_{k+1} := \tomega_k/\sqrt{\theta_\omega}$
				\State $\tblambda_k := \bO$
			\Else
				\State $\blambda_{k+1} := \blambda_k$\,,\quad $\tomega_{k+1} = \theta_\omega \cdot \tomega_k$
			\EndIf
		\EndIf
	\EndFor
\end{algorithmic}
The steps of the algorithm are clear since it is -- with the only exception in line 4 -- fully identical to the algorithm of ALM discussed in Section~\ref{Sec:ALM}.

We briefly discuss the termination condition in line 6. We find that \eqref{eqn:Newton_residual} coincides with the KKT condition of \eqref{eqn:PenaltyProblem} for $\omega>0$ and with the KKT conditions of \eqref{eqn:ConstrainedProblem} when $\omega=0$. The only differences appear in terms of $\tblambda$. Thus, when $\tblambda$ has converged to zero, which it does as formerly discussed, then the KKT system for the respective value of $\omega$ is solved up to the residual tolerance $\tol>0$.

\section{Outlook: Sophisticated iterations for the penalty program}
We discussed three methods: Penalty method, ALM, and modified ALM. All these methods employ a Newton iteration to solve penalty problems of some form
\begin{align*}
	\min_{\by \in \R^p}\quad g(\by) + \frac{1}{2\cdot\tomega} \cdot \|d(\by)\|_2^2\,.\tageq\label{eqn:DummyPenaltyProblem}
\end{align*}
The exact functions for $g,d$ and the exact iterates for $\by$ are method-dependent and shall not be further specified for the purpose of this section.

For sake of simplicity, in the former sections we always used a line-seearch Newton method for solving the penalty program. However, practical methods may use enhanced iterations for solving the penalty problem. These sophistications are: An adaptive choice of the penalty parameter $\tomega$; the use of second-order corrections within the line search.

\paragraph{Adaptive choice of the penalty parameter}
The value of $\tomega$ must be just of the right order of magnitude in order to achieve the following two antagonistic goals:
\begin{itemize}
	\item $\tomega$ must be sufficiently large so that \eqref{eqn:DummyPenaltyProblem} is well-scaled, in order to admit large steps within the Newton iteration.
	\item $\tomega$ must be sufficiently small to yield boundedness of \eqref{eqn:DummyPenaltyProblem}, otherwise a local minimizer may not exist.
\end{itemize}
In order for the outer iteration $k$ to converge, there is usually an upper bound for $\tomega$. Let us call this bound $\hat{\omega}_k$. Choosing $\tomega:=\hat{\omega}_k$ may not be sufficient to guarantee boundedness of \eqref{eqn:DummyPenaltyProblem}. Thus, it is of benefit for a practical implementation that $\tomega$ can be decreased within the inner iteration. Also, when $\tomega< \hat{\omega}_k$ and the line search makes slow progress, we may think of increasing $\tomega$ by using a suitable adaptive scheme, e.g., as presented in \cite{Amudson}. While it is of interest for our paper that we pointed out there is freedom in adapting $\tomega$, a further discussion of how this adaption is performed is beyond the scope of this paper.

\paragraph{Second-order corrections within the line search}

Second-order corrections (SOC) are often motivated as a way for avoiding the Maratos effect \cite{NumOpt}. Throughout this paper we did only consider unconstrained minimization techniques \cite{SUMT} that do not suffer from the Maratos effect by construction. However, it is acknowledged that in practice the line search benefits from second-order corrections \cite{NumOpt}.

We spend care that the SOC is formulated in a correct way with respect to our problem. Usually, consistency is achieved by the fact that the Newton iteration drives $c(\bx)= \bO$. But in our case this equation will never be satisfied, making it sophisticated to construct a consistent SOC. Below we give a consistent SOC for the modified ALM.

\largeparbreak

For the modified ALM we applied a quadratic penalty to the constraint
$$\tilde{c}(\by)\equiv c(\bx) + \omega \cdot (\blambda_k + \tblambda) = \bO\,.$$
This equation has a solution but the penalty method will not converge to it. Instead, as we find from \eqref{eqn:Newton_residual}, the iteration will eventually solve
\begin{align}
c(\bx) + \omega \cdot \blambda_k + (\omega + \tomega_k) \cdot \tblambda = \bO\,. \label{eqn:17}
\end{align}
Therefore, at given a trial point $\bw(\alpha) := (\bx+\alpha \cdot \Delta\bx,\tblambda + \alpha \cdot \Delta\tblambda)$, it makes sense to consider its second-order corrected trial point
\begin{align*}
	\bw(\alpha)^{\text{SOC}} := \bw(\alpha) + (\Delta\bx^{\text{SOC}},\Delta\tblambda^{\text{SOC}})\,,
\end{align*}
that we construct from the solution of the following linear system:
\begin{align*}
	\begin{bmatrix}
		\bS & \nabla c(\bx)\\
		\nabla c(\bx)\t & - (\omega+\tomega) \cdot \bI
	\end{bmatrix} \cdot \begin{pmatrix}
		\Delta\bx^{\text{SOC}}\\
		-\Delta\tblambda^{\text{SOC}}
	\end{pmatrix} = - \begin{pmatrix}
		\bO\\
		c\big(\bx(\alpha)\big) + \omega \cdot \blambda_k + (\omega+\tomega) \cdot \tblambda(\alpha)
	\end{pmatrix}
\end{align*}
In \cite{IPOPT} the authors propose to use multiple SOC steps. That is, a point $\bw^{\text{SOC}}$ is considered as trial point $\bw(\alpha)$, and is again incremented with a correction step following the above construction.

The matrix $\bS \in \R^{n \times n}$ must be symmetric positive definite and should represent the problem's scaling. The matrix determines the induced $\bS$-norm in which $\Delta\bx^{\text{SOC}}$ is minimized subject to solving a linearization of \eqref{eqn:17}.

\section{Conclusions}
The modified ALM is an iterative method that can solve programming problems with equality constraints and problems with large quadratic penalty terms. The modified ALM is a generalization of the ALM and inherits all convergence properties of the ALM. Like the ALM, the modified ALM solves a sequence of subproblems whose penalty-parameters remain of moderate size under the same suitable conditions as for the ALM. This is advantageous for the numerical solution of the Newton system and for achieving large step sizes in the line search.

Future work is related to incorporating this method into SQP methods and Interior-Point methods for the general non-linear programming problem.

\FloatBarrier

\bibliography{MALM_bib}
\bibliographystyle{plain}

\end{document}